\documentstyle[amsthm,amsfonts,amssymb,amsgen,amsmath,amsopn,verbatim,hhline,tabularx,11pt]{article}

\theoremstyle{plain}
\newtheorem{thm}{Theorem}[section]
\newtheorem{lem}[thm]{Lemma}

\newtheorem{cor}[thm]{Corollary}
\newtheorem{conj}[thm]{Conjecture}

\theoremstyle{definition}

\newcommand{\lra}{\longrightarrow}
\newcommand{\wg}{\wedge}
\newcommand{\bwg}{\bigwedge}
\newcommand{\llra}{\longleftrightarrow}
\newcommand{\tZ}{{\tilde{Z}}}
\newcommand{\Z}{{\mathbb Z}}
\renewcommand{\P}{{\mathbb P}}
\newcommand{\Q}{{\mathbb Q}}

\DeclareMathOperator{\sgn}{{sgn}}
\DeclareMathOperator{\cyc}{{cyc}}

\newcommand{\eps}{{\varepsilon}}
\newcommand{\pa}{{\partial}}
\newcommand{\ga}{{\alpha}}

\newcommand{\gb}{{\beta}}
\newcommand{\gd}{{\delta}}
\newcommand{\gs}{{\sigma}}

\newcommand{\os}{\overset}
\newcommand{\ot}{\otimes}
\newcommand{\bul}{{\bullet}}

\newcommand{\caB}{{\mathcal B}}
\newcommand{\CH}{{\mathcal CH}}

\newcommand{\caG}{{\mathcal G}}
\newcommand{\caT}{{\mathcal T}}

\newcommand{\lms}{\longmapsto}

\begin{document}

\title{Goncharov's Relations in Bloch's higher Chow Group $CH^3(F,5)$
}
\author{Jianqiang Zhao\footnote{Partially supported by NSF grant DMS0139813}}
\date{}
\maketitle
\begin{center}
Department of Mathematics, Eckerd College, St. Petersburg, FL
33711
\end{center}

\noindent{\bf Abstract.} In this paper we will prove Goncharov's
22-term relations in the linearized version of the Bloch's higher
Chow group $CH^3(F,5)$ using linear fractional cycles of
Bloch-Kriz-Totaro under the Beilinson-Soul\'es
vanishing conjecture that $CH^2(F,n)=0$ for $n\ge 4$.

\section{Introduction}
Around 1980 Goncharov defined his polylogarithmic
(cohomological) motivic complex over an arbitrary field $F$:
$$\Gamma(F,n): \caG_n(F) \os{\gd_n}{\lra} \caG_{n-1}(F) \ot F^\times
\os{\gd_{n-1}}{\lra} \cdots
\os{\gd_3}{\lra} \caG_2(F) \ot \bwg^{n-2} F^\times \os{\gd_2}{\lra}
 \bwg^n F^\times,$$
where $\caG_n$, denoted by $\caB_n$ by Goncharov, is placed at
degree 1. To save space we here only point out that $\caG_n(F)$
are quotient groups of $\Z[\P^1_F]$ and refer the interested readers
to \cite[p.49]{Ggalois} for the detailed definition of these groups.

On the other hand, currently there are two versions of higher Chow
groups available: simplicial and cubical, which are isomorphic
(cf.~\cite{Lev}). We will recall the cubical version in \S 2 and
use it throughout this paper.

Define the $\Z$-linear map $\gb_2: \Z[\P^1_F]\lra \bwg^2 F^\times$
by $\gb_2(\{x\})=(1-x)\wg x$ for $x\ne 0,1$ and $\gb_2(\{x\})=0$
for $x=0,1$. Let $B_2(F)$ be the Bloch group defined as the
quotient group of  $\ker(\gb_2)$ by the 5-term relations satisfied
by the dilogarithm. In \cite{GM} Gangl and M\"uller-Stach prove
that there is a well-defined map to the higher Chow group
$$\bar\rho_2:B_2(F)_\Q \lra CH^2(F,3)_\Q,$$
where we denote $G_\Q=G\otimes\Q$ for any abelian group $G$. The
essential difficulty of the proof lies in showing that
$\rho_2:\Z[\P_F^1]\lra CH^2(F,3)$ sends 5-term relations to 0,
where for $a\in F$ we assign $\rho_2(\{a\})$ the linear fractional
cycle $C_a^{(2)}$ of Totaro \cite{T}, generalized by Bloch and
Kriz \cite{BK}. For $m\ge 2$ these cycles are defined as
$$C_a^{(m)}=\Bigl[x_1,\dots,x_{m-1},1-x_1,
1-\frac{x_2}{x_1},\dots,1-\frac{x_{m-1}}{x_{m-2}},
1-\frac{a}{x_{m-1}}\Bigr]\in CH^m(F,2m-1).$$
It is believed that $\bar\rho_2$ gives rise to an isomorphism
because of the following results of Suslin (cf. \cite{Sus1,Sus2}):
$B_2(F)_\Q \cong K_3^{\text{ind}}(F)_\Q \cong CH^2(F,3)_\Q$
where $K_3^{\text{ind}}(F)$ is the indecomposable part of $K_3(F)$.

One naturally hopes to carry the above to the higher Chow groups
$CH^m(F,2m-1)_\Q$ for $m\ge 3$. One can start by defining
$\caB_m(F)$ as the the subgroup $\ker \gd_m$ of $\caG_m$ because
it is known that $\caB_2(F)\cong B_2(F)$ for number fields $F$, at
least modulo torsion. (There are some other ways to define these
groups, see \cite{Z}). One then has
\begin{conj} For $m\ge 3$,
$$\caB_m(F)_\Q\cong CH^m(F,2m-1)_\Q.$$
\end{conj}
Even for $m=3$ the current state of knowledge requires
modifications of the groups on both sides. For example, we do not
yet have a very good understanding of the relation group of
$\caG_3(F)$ although we expect it is equal to $R_3(F)$ which is
generated by the following relations:

(1) $\{x\}-\{x^{-1}\}$, $x\in F$;

(2) $\{x\}+\{1-x\}+\{1-x^{-1}\}-\{1\}$, $x\in F$;

(3) Goncharov's 22-term relations: for any $a,b,c\in \P^1_F$
\begin{multline*}
R(a,b,c)=\{-abc\}+\bigoplus_{\cyc(a,b,c)}\Bigl(
\{ca-a+1\}+\Bigl\{\frac{ca-a+1}{ca}\Bigr\}
-\Bigl\{\frac{ca-a+1}{c}\Bigr\} \\
+\Bigl\{\frac{a(bc-c+1)}{-(ca-a+1)}\Bigr\}
+\Bigl\{\frac{bc-c+1}{b(ca-a+1)}\Bigr\}
+\{c\}-\Bigl\{\frac{bc-c+1}{bc(ca-a+1)}\Bigr\}-\{1\}\Bigr).
\end{multline*}
We thus define $B_3(F)$ as $\ker(\gb_3: \Z[\P^1_F]\lra
B_2(F)\otimes F^\times)/R_3(F)$. This is well defined by a result
of Goncharov \cite{Ggalois}. We then replace the group
$CH^3(F,5)_\Q$ by $\CH^3(F,5)_\Q$ (see \S~\ref{setup}) which is
isomorphic to $\CH^3(F,\bul)$ by some mild conjecture. The
ultimate goal of our work is to prove
\begin{conj} Let $F$ be a field. Then
$$B_3(F)_\Q\cong \CH^3(F,5)_\Q\cong CH^3(F,5)_\Q. $$
\end{conj}
Define the map
$$\aligned
\rho_3:\Q[\P_F^1] &\lra \CH^3(F,5)_\Q \\
    \{a\} &\lms C_a^{(3)}.
\endaligned$$
Let $T(a)=\{a\}+\{1-a\}+\{1-a^{-1}\}$. By \cite[Thm.~2.9(b)]{GM}
we know that for any $a,b\ne 0,1$ in $F$ we have
$\rho_3(T(a))=\rho_3(T(b))$. We denote this cycle by $\eta$. The
main purpose of this paper is to show that if we replace $\{1\}$
by $\eta$ in relation (3) then it is sent to 0 under $\rho_3$ when
none of the terms is $\{0\}$ or $\{1\}$. Note that Gangl and
M\"uller-Stach has done the same for (1) and they even prove the
Kummer-Spence relations which are special cases of (3). Naturally,
our work builds on theirs. The proof of relation (2) in
$\CH^3(F,5)_\Q$ is still open as of now.

To simplify exposition we disregard torsions throughout this paper.
In fact, all the results are still valid if we modulo 4-torsions only.

I was attracted to \cite{GM} by a few very helpful conversations
with Owen Patashnick. I gratefully acknowledge the hospitality
provided by the Mathematics Department of Duke University during
my visit. I also want to thank Dick Hain for his interest, help
and encouragement and Herbert Gangl for pointing out a serious
mistake about the admissibility of cycles in the first draft of
this paper.

\section{The setup}\label{setup}
Let $F$ be an arbitrary field.  The algebraic $n$-cube
$$\square^n=(\P_F^1\setminus\{1\})^n$$
has $2^n$ codimension one faces given by $\{t_i=0\}$ and $\{t_i=\infty\}
$
for $1\le i\le n$. We have the boundary map
$$ \pa =\sum_{i=1}^n (-1)^{i-1} (\pa_i^0- \pa_i^\infty)$$
where $\pa_i^a$ denotes the restriction map on face $t_i=a$.
Recall that for a field $F$ one let $Z^p_c(F,n)$ (subscript $c$
for ``cubical'') be the quotient of the group of admissible
codimension $p$ cycles in $\square^n$ by the subgroup of
degenerate cycles as defined in \cite[p.180]{T}. {\em Admissible}
means that the cycles have to intersect all the faces of any
dimension properly. Levine \cite{Lev} shows that the $n$-th
homology group of the resulting complex $Z^p_c(F, \bullet)$ is
isomorphic to Bloch's higher Chow group $CH^p(F,n)$. This
establishes the isomorphism between the cubical and simplicial
version of Bloch's higher Chow groups. Furthermore, Bloch
\cite{Bl} constructs a rational alternating version $C^p(F,n)$ of
$Z^p_c(F, \bullet)$ (see also \cite[\S2]{GM}) whose homological
complex
$$C^m(F,\bul): \
\cdots\lra C^m(F,2m)\lra C^m(F,2m-1)\lra \cdots \lra C^m(F,m)\lra 0$$
still computes $CH^p(X,n)_\Q$ as proved in \cite{Lev}. The properties
of the elements in $C^m(F,n)$ are essentially
encoded in the following equation: for any choice of
$\gd_1,\dots,\gd_n=\pm 1$ and any permutation $\gs$ of $\{1,\dots,n\}$
$$\big[f_1^{\gd_1},\dots, f_n^{\gd_n}\big]=
 \sgn(\gs)\Bigl(\prod_{i=1}^n \gd_i \Bigr)
    \big[f_{\gs(1)},\dots, f_{\gs(n)}\big].$$

To simplify computation Gangl and M\"uller-Stach further modify
the complex $C^3(F,\bul)$ by taking the quotient by an acyclic
subcomplex $S^3(F,\bul)$. (See their paper for the definition.
Also note that acyclicity is proved under Beilinson-Soul\'es
conjecture $CH^2(F,n)=0$ for $n\ge 4$). Following them we call
cycles in $S^3(F,\bul)$ {\em negligible} and denote the quotient
complex by $A^3(F,\bul)$. We further put
$\CH^3(F,n)=H_n(A^3(F,\bul))$ (note the different fonts). Hence
$$\CH^3(F,n)\cong CH^3(F,n)$$
under the conjecture $CH^2(F,n)=0$ for $n\ge 4$.

\section{Some lemmas}
We will mostly follow the notation system in \cite{GM} except
that we denote
$$\{a\}_c=\Bigl[x,y,1-x,1-\frac{y}{x},1-\frac{a}{y}\Bigr].$$
This is denoted by $C_a$ in \cite{GM}. The subscript $c$ here is
for ``cubical''.

\begin{lem}\label{fgh}{\em (Gangl-M\"uller-Stach)}
Let $f_i$ ($i=1,2,3,5$) be rational functions and $f_4(x,y)$
be a product of fractional linear transformations of the form
$(a_1x+b_1y+c_1)/(a_2x+b_2y+c_2)$.  We assume that all the cycles
in the lemma are {\em admissible} and write
$$Z(f_1,f_2)=[f_1,f_2,f_3,f_4,f_5]
=[f_1(x),f_2(y),f_3(x),f_4(x,y),f_5(y)]$$
if no confusion arises.

\noindent(i) If $f_4(x,y)=g(x,y)h(x,y)$ then
$$[f_1,f_2,f_3,f_4,f_5]=[f_1,f_2,f_3,g,f_5]+[f_1,f_2,f_3,h,f_5].$$

\noindent(ii) Assume that $f_1=f_2$ and
that for each non-constant solution $y=r(x)$ of $f_4(x,y)=0$
and $1/f_4(x,y)=0$ one has $f_2(r(x))=f_2(x)$.\begin{quote}
(a) If $f_3(x)=g(x)h(x)$ then
$$[f_1,f_2,f_3,f_4,f_5]=[f_1,f_2,g,f_4,f_5]+[f_1,f_2,h,f_4,f_5].$$

(b) Similarly, if $f_5(y)=g(y)h(y)$ then
$$[f_1,f_2,f_3,f_4,f_5]=[f_1,f_2,f_3,f_4,g]+[f_1,f_2,f_3,f_4,h].$$

(c) If $f_1=f_2=gh$ and
$g(r(x))=g(x)$ or $g(r(x))=h(x)$ then
\begin{equation}\label{1steq}
2 Z(f_1,f_2)=Z(g,f_2)+Z(h,f_2)+Z(f_1,g)+Z(f_1,h)
\end{equation}
and
\begin{equation}\label{2ndeq}
Z(f_1,f_2)=Z(g,g)+Z(h,h)+Z(h,g)+Z(g,h).
\end{equation}
\end{quote}
\end{lem}
\begin{proof}
(i) is  contained in Lemma 2.8(b) of \cite{GM}.
(ii) can be proved using the same idea as in the
proof of \cite[Lemma 2.8(a)(c)]{GM}.
\end{proof}

\begin{lem}\label{keylem}  Assume that $f_i$, $i=1,2,3,5$, are
rational functions of one variable and $p_4$ and $q_4$ are
rational functions of two variables.
Assume that the only non-constant solution of
$p_4(x,y)=0,\infty$ is $y=x$ and the same for $q_4(x,y)$.

(i) If $f_3=gh$ then
$$\aligned
\ [f_1,f_2,f_3, p_4, f_5] + [f_2,f_1,f_3, q_4, f_5]
=&[f_1,f_2,g, p_4, f_5]+ [f_2,f_1, g, q_4, f_5]\\
+&[f_1,f_2,h, p_4, f_5]+ [f_2,f_1, h, q_4, f_5]
\endaligned$$
if all cycles are admissible. A similar result holds if $f_5=gh$.

(ii) If $f_2=gh$ then
$$\aligned
\ [f_1,f_2,f_3, p_4, f_5] + [f_2,f_1,f_3, q_4, f_5]
=&[f_1,g, f_3,p_4, f_5]+ [g,f_1, f_3, q_4, f_5]\\
+&[f_1,h, f_3,p_4, f_5]+ [h,f_1, f_3, q_4, f_5]
\endaligned$$
if all cycles are admissible.
\end{lem}
\begin{proof} (i) Write $[F_1,\dots,F_6]=[F_1(x),F_2(y),\dots,F_5(x,y),F_6(y)]$
and let
$$W=\Bigl[f_1,f_2,\frac{z-g(x)h(x)}{z-g(x)},z,p_4, f_5\Bigr]+
\Bigl[f_2,f_1,\frac{z-g(x)h(x)}{z-g(x)},z,q_4, f_5\Bigr].$$
Taking the boundary we get the desired result because the cycle
$$V=\Bigl[f_1(x),f_2(x),\frac{z-g(x)h(x)}{z-g(x)},z, f_5(x)\Bigr]$$
cancels with
$$-V=\Bigl[f_2(x),f_1(x),\frac{z-g(x)h(x)}{z-g(x)},z,f_5(x)\Bigr]$$
by skew-symmetry.

(ii) This is similar to (i) if we set
$$W=\Bigl[f_1(x), \frac{z-g(y)h(y)}{z-g(y)},z,f_3, p_4, f_5\Bigr]
-\Bigl[\frac{z-g(x)h(x)}{z-g(x)},z,f_1(y),f_3,q_4, f_5\Bigr].$$
\end{proof}

\begin{cor}\label{key}
If the conditions in the lemma are all satisfied then for
$\ga\in F$
$$[f_1,f_2,\ga f_3, p_4, f_5] + [f_2,f_1,\ga f_3, q_4, f_5]
=[f_1,f_2,f_3, p_4, f_5] + [f_2,f_1,f_3, q_4, f_5].$$
A similar result holds if the the constant $\ga$ is in front
of $f_5$.
\end{cor}

The next computational lemma is easy. Here and in
what follows we formally extend the definition of $\{?\}_c$ to
include $\{0\}_c=\{\infty\}_c=0$.
\begin{lem}\label{stuv} For all $s,t,u,v\in F$ the admissible cycle
$$\Bigl[x,y, \frac{1-sx}{1-tx}, 1-\frac{y}{x}, \frac{u-y}{v-y}\Bigr]
=\{us\}_c-\{vs\}_c-\{ut\}_c+\{vt\}_c.$$
Similarly, the admissible cycle
$$\Bigl[x,y, \frac{s-x}{t-x}, 1-\frac{x}{y}, \frac{1-uy}{1-vy}\Bigr]
=\{us\}_c-\{vs\}_c-\{ut\}_c+\{vt\}_c.$$
\end{lem}
\begin{proof}
By Lemma \ref{fgh}(ii) we only need to show
\begin{equation}\label{us}
\Bigl[x,y, 1-sx, 1-\frac{y}{x},1-\frac{u}{y}\Bigr]=\{us\}_c
\end{equation}
which follows easily from a substitution $(x,y)\mapsto (x/s, y/s)$
if $s\ne 0$. If $s=0$ then \eqref{us} is trivial.
The second equation follows from the obvious
substitution $(x,y)\mapsto (y,x)$.
\end{proof}

\section{Goncharov's relations} \label{mainsec}
Let $\caT(a)=\{a\}_c+\{1-a\}_c+\{1-a^{-1}\}_c$.
By \cite[Thm.~2.9(b)]{GM} we know that for any $a,b\ne 0,1$ in $F$
we have $\caT(a)=\caT(b)$. We denote this cycle by $\eta$.

\begin{thm}\label{main}
Goncharov's 22 term relations hold in $\CH^3(F,5)$: for any $a,b,c\in
\P_F^1$
\begin{multline}\label{Rabc=0}
R(a,b,c)=\{-abc\}+\bigoplus_{\cyc(a,b,c)}\Bigl(
\{ca-a+1\}+\Bigl\{\frac{ca-a+1}{ca}\Bigr\}
-\Bigl\{\frac{ca-a+1}{c}\Bigr\} \\
+\Bigl\{\frac{a(bc-c+1)}{-(ca-a+1)}\Bigr\}
+\Bigl\{\frac{bc-c+1}{b(ca-a+1)}\Bigr\}
+\{c\}-\Bigl\{\frac{bc-c+1}{bc(ca-a+1)}\Bigr\}-\eta\Bigr)=0,
\end{multline}
where $\cyc(a,b,c)$ means cyclic permutations of $a,b$ and $c$,
provided that none of terms is $\{0\}$ or $\{1\}$ except for
$\eta$ (non-degeneracy condition). Here we drop the subscript $c$
for the cycle notation $\{ ?\}_c.$
\end{thm}
\begin{proof}
To make the proof explicit we will carry it out in a series of
steps. Throughout the proof we will use $\{1/t\}=\{t\}$ repeatedly
without stating it explicitly. As Gangl pointed out to the author
the major difficulty is to guarantee that all the cycles we use
lie in the ``admissible world''. Due to its length and pure
computational feature we put the proof of admissibility of all the
cycles appearing in this paper in the online supplement
\cite{internet} except for one cycle in Step (2) where we spell
out all the details to provide the readers the procedure how we do
the checking in general.

\medskip
\noindent
{\bf Step} (1). Construction of $\{k(c)\}$.

\medskip
Let $f(x)=x$, $A(x)=(ax-a+1)/a$ and $B(x)=bx-x+1$ . Let
$k(x)=B(x)/abxA(x)$ and $l(y)=1-(k(c)/k(y))$. Then taking
$\mu=-(ab-b+1)/a$ we can write
\begin{equation}\label{rmab}
\{k(c)\}=\Bigl[x,y,1-x,1- \frac{y}{x},1- \frac{k(c)}{y}\Bigr]=
 \Bigl[\frac{abx}{\mu}, \frac{aby}{\mu} ,1-x,1- \frac{y}{x},1- \frac{k(c)}{y}\Bigr]
\end{equation}
by Lemma \ref{fgh}(ii) because all of the following cycles are
admissible and negligible
\begin{align*}
 \Bigl[\frac{ab}{\mu},y,f_3,f_4,f_5\Bigr],\quad
 \Bigl[x,\frac{ab}{\mu},f_3,f_4,f_5\Bigr],\quad
 \Bigl[\frac{ab}{\mu},\frac{ab}{\mu},f_3,f_4,f_5\Bigr],
\end{align*}
where $(f_3,f_4,f_5)=(1-x,1-y/x,1-k(c)/y)$. Here for the last
cycle we need to use the fact that
\begin{equation}\label{kcne1}
1-k(c)=\frac{(c-1)(1+abc)}{abcA(c)}\ne 0.
\end{equation}

Next by using the transformation $(x,y)\mapsto (k(x),k(y))$ we get
$$4\{k(c)\}=\Bigl[\frac{B(x)}{\mu xA(x)},\frac{B(y)}{\mu yA(y)},
 1-k(x),1-\frac{k(y)}{k(x)},l(y)\Bigr]
 =Z\Bigl(\frac{B}{\mu fA},\frac{B}{\mu fA}\Bigr).$$
Here for any two rational functions $f_1$ and $f_2$ of one
variable we set
$$Z(f_1,f_2)=\Bigl[f_1(x),f_2(y),1-k(x),1-\frac{k(y)}{k(x)}, l(y)
\Bigr].$$

\medskip
\noindent
{\bf Step} (2). The key reparametrization and a
simple expression of $\{k(c)\}$.

\medskip

We first observe that under the involution $x \os{\rho_x}{\llra}
-A(x)/B(x)$ we have
$$\aligned
k(x)\os{\rho_x}{\llra} & k(x), \quad \frac{x-1}{x}
\os{\rho_x}{\llra} \frac{abx+1}{aA(x)},\quad
1-\frac{\mu x}{A(y)B(x)}\os{\rho_x}{\llra} \frac{y-x}{A(y)}\\
B(x)\os{\rho_x}{\llra} &\frac{-\mu}{B(x)}, \quad \frac{B(x)}{x}
\os{\rho_x}{\llra} \frac{\mu}{A(x)}, \frac{A(x)}{x}
\os{\rho_x}{\llra}  \frac{-\mu x}{A(x)}.
\endaligned$$
Next if applying both $\rho_x$ and $\rho_y$ (denoted by
$\rho_{x,y}$) then we get
\begin{equation}\label{rhoxys}
\begin{split}
1-\frac{x}{y} \os{\rho_{x,y}}{\llra} \frac{\mu(x-y)}{A(y)B(x)},
&\quad
\frac{y-x}{yB(x)} \os{\rho_{x,y}}{\llra} \frac{y-x}{A(y)},\\
\frac{A(y)}{y}\Bigl(1-\frac{\mu x}{A(y)B(x)}\Bigr)
\os{\rho_{x,y}}{\llra} & B(x)\Bigl(1-\frac{\mu x}{A(y)B(x)}\Bigr).
\end{split}
\end{equation}
By Lemma \ref{fgh}(ii)
\begin{multline}\label{AA}
4\{k(c)\} = Z\Bigl(\frac{\mu fA}{B}, \frac{\mu fA}{B}\Bigr)
    =Z( A,A)+Z\Bigl(\frac{\mu f}{B},A\Bigr)
    +Z\Bigl(A,\frac{\mu f}{B}\Bigr)
    +Z\Bigl(\frac {\mu f}{B},\frac{\mu f}{B}\Bigr) \\
= Z(A,A)+\rho_x Z(A,A)+\rho_y Z(A,A)+\rho_{x,y}Z(A,A)=4Z(A,A).
\end{multline}
We end this step by showing that $Z_A=Z(A,A)$ is admissible. Note
that
\begin{align} \label{1-kx}
1-k(x)=&\frac{(x-1)(1+abx)}{abxA(x)},\\
1-\frac{k(y)}{k(x)}=&\frac{(y-x)(yB(x)+A(x))}{yA(y)B(x)}=\frac{(y-x)(xB(y)+A(y))}{yA(y)B(x)}.
\label{1-kyx}
\end{align}
We have
\begin{align*}
    \pa_1^0(Z_A)&\subset \{t_4=1\},\quad
   \pa_1^\infty(Z_A)\subset \{t_3=1\},\\
   \pa_2^0(Z_A)&\subset \{t_5=1\},\quad
   \pa_2^\infty(Z_A)\subset \{t_4=1\},\quad
   \pa_3^\infty(Z_A)\subset \{t_4=1\},\\
   \pa_4^\infty(Z_A)&\subset \{t_3=1\}\cup \{t_5=1\},\quad
   \pa_5^\infty(Z_A)\subset \{t_4=1\},\\
    \pa_3^0(Z_A)&=\Bigl[\frac{1}{a},A(y),1-k(y),l(y)\Bigr]
    +\Bigl[A\Big(\frac{-1}{ab}\Big),A(y),1-k(y),l(y)\Bigr],\\
   \pa_4^0(Z_A)&=\Bigl[A(y),A(y),1-k(y),l(y)\Bigr]
   +\Bigl[\frac{\mu y}{B(y)},A(y),1-k(y),l(y)\Bigr],\\
   \pa_5^0(Z_A)&=\Bigl[A(x),A(c),1-k(x),l(x)\Bigr]
   +\Bigl[A(x),A(y_2),1-k(x),l(x)\Bigr],
\end{align*}
where the last equation comes from the two solutions of $l(y)=0$:
\begin{equation*}
y_1=c \quad\text{ and }  \quad
y_2=-\frac{ac-a+1}{a(bc-c+1)}=-\frac{A(c)}{B(c)}=\rho_c(c).
\end{equation*}
By non-degeneracy assumption and
\begin{align}
A(y_2)= &\rho_c(A(c))=c\mu /B(c)\ne 0,\infty, \nonumber \\
B(y_2)= &\rho_c(B(c))=-\mu/B(c) \ne 0,\infty, \label{By2}
\end{align}
it suffices to show the following cycles are admissible:
\begin{alignat*}{3}
    L&:=\Bigl[A(y),1-k(y),l(y)\Bigr],\quad&
    L'&:=\Bigl[A(y),A(y),1-k(y),l(y)\Bigr],\\
    L''&:=\Bigl[\frac{\mu y}{B(y)},A(y),1-k(y),l(y)\Bigr].
\end{alignat*}

\begin{itemize}
\item $L$ is admissible. Because $l(y)=1-yB(c)A(y)/cA(c)B(y)$ we
have
$$\pa_1^0(L)\subset \{t_3=1\},\quad
 \pa_1^\infty(L)\subset\{t_2=1\},\quad
 \pa_2^\infty(L)\subset\{t_3=1\},\quad
 \pa_3^\infty(L)\subset\{t_2=1\}.
$$
Moreover, by non-degeneracy assumption we see that by
\eqref{kcne1} and \eqref{1-kx}
\begin{align*}
A(1)=&\frac{1}{a}\ne 0, \qquad k(1)=1, \qquad  l(1)=1-k(c)\ne 0,\\
 A\Big(\frac{-1}{ab}\Big)=&\frac{\mu}{b}\ne 0,\quad k(-1/ab)=1,
 \quad l\Big(\frac{-1}{ab}\Big)= 1-k(c)\ne 0,\\
 aby_2+1&=\frac{(1-c)(ab-b+1)}{bc-c+1} \ne 0.
\end{align*}
Thus both $\pa_2^0(L)=[A(1),l(1)]+[A(-1/ab),l(-1/ab)]$ and
$\pa_3^0(L)=[A(c),1-k(c)]+[A(y_2),1-k(c)]$ are clearly admissible
by non-degeneracy assumption.

\item $L'$ is admissible. This follows from the above proof for
$L$.

\item $L''$ is admissible. This also follows from the proof for
$L$ because $\mu y/B(y)\ne 0,\infty$ when $y=1,-1/ab,c, y_2$ by
\eqref{By2}.

\end{itemize}

\medskip
\noindent
{\bf Step} (3). Some admissible cycles for
decomposition of $\{k(c)\}$.

\medskip

In order to decompose $Z(A,A)$ we define the following
admissible cycles
$$\aligned
Z_1(A,A)=&\Bigl[\frac{(b-1)A(x)}{\mu},\frac{(b-1)A(y)}{\mu},
        \frac{x-1}{x}, \frac{y-x}{A(y)}, l(y)\Bigr],\\
Z_2(A,A)=&\Bigl[\frac{(b-1)A(x)}{\mu},\frac{(b-1)A(y)}{\mu},
        \frac{x-1}{x}, \Bigl(\frac{A(y)}{y}\Bigr)
    \Bigl(1- \frac{\mu x}{A(y)B(x)}\Bigr),l(y)\Bigr],\\
Z_3(A,A)=&\Bigl[\frac{(b-1)A(x)}{\mu},\frac{(b-1)A(y)}{\mu},
        \frac{abx+1}{abA(x)}, \frac{y-x}{A(y)},l(y)\Bigr],\\
Z_4(A,A)=&\Bigl[\frac{(b-1)A(x)}{\mu},\frac{(b-1)A(y)}{\mu},
        \frac{abx+1}{abA(x)}, \Bigl(\frac{A(y)}{y}\Bigr)
        \Bigl(1- \frac{\mu x}{A(y)B(x)}\Bigr),l(y)\Bigr].
\endaligned$$
We now use Lemma~\ref{fgh}(ii(c)) to remove the coefficients in
front of $A(x)$ and $A(y)$ in $Z_1(A,A)$ and $Z_3(A,A)$,
Lemma~\ref{fgh}(i) to remove the factor $A(y)/y$ from the fourth
coordinate of $Z_2(A,A)$, and Lemma~\ref{fgh}(ii(a)) to remove the
coefficient $1/b$ in front of the third coordinates of $Z_4(A,A)$,
and finally get:
$$\aligned
Z_1(A,A)=&\Bigl[A(x),A(y),\frac{x-1}{x}, \frac{y-x}{A(y)}, l(y)\Bigr],\\
Z_2(A,A)=&\Bigl[\frac{(b-1)A(x)}{\mu},\frac{(b-1)A(y)}{\mu},
        \frac{x-1}{x}, 1- \frac{\mu x}{A(y)B(x)},l(y)\Bigr],\\
Z_3(A,A)=&\Bigl[A(x),A(y),\frac{abx+1}{abA(x)},
        \frac{y-x}{A(y)}, l(y)\Bigr],\\
Z_4(A,A)=&\Bigl[\frac{(b-1)A(x)}{\mu},\frac{(b-1)A(y)}{\mu},
        \frac{abx+1}{aA(x)}, \Bigl(\frac{A(y)}{y}\Bigr)
        \Bigl(1- \frac{\mu x}{A(y)B(x)}\Bigr),l(y)\Bigr].
\endaligned$$
It is not too hard to verify that all the cycles appearing in the
above are admissible.

Now we can break up the fourth coordinate of $Z(A,A)$ according to
Lemma~\ref{fgh}(i) and get
\begin{multline*}
Z(A,A)=Z'(A,A)+Z''(A,A)=
\Bigl[\frac{(b-1)A(x)}{\mu},\frac{(b-1)A(y)}{\mu},
        1-k(x), \frac{y-x}{A(y)},l(y)\Bigr] \\
+\Bigl[\frac{(b-1)A(x)}{\mu},\frac{(b-1)A(y)}{\mu}, 1-k(x),
\frac{A(y)}{y}\Bigl(1- \frac{\mu x}{A(y)B(x)}\Bigr),l(y)\Bigr].
\end{multline*}
Also by Lemma \ref{fgh}(ii) we find that
$$Z'(A,A)=Z_1(A,A)+Z_3(A,A).$$
However, the conditions in Lemma~\ref{fgh}(i) are not all
satisfied by $Z''(A,A)$. To decompose the third coordinates of
$Z''(A,A)$ we combine
$$Z''(A,A)=\rho_y Z''(A, A)=
\Bigl[\frac{(b-1)A(x)}{\mu},\frac{(b-1)y}{B(y)}, 1-k(x),
\frac{\mu(x-y)}{A(y)B(x)},l(y)\Bigr]$$ and
$$Z''(A,A)=\rho_x Z''( A, A)=
\Bigl[\frac{(b-1)x}{B(x)},\frac{(b-1)A(y)}{\mu}, 1-k(x),
1-\frac{x}{y},l(y)\Bigr]$$ and use Lemma~\ref{keylem}(i) to get
$$2 Z''(A,A) =(\rho_x+\rho_y)\Big(Z_2(A,A)+Z_4(A,A)\Big)
    =2Z_2(A,A)+2Z_4(A,A).$$
Here the properties of substitutions $\rho_x$ and $\rho_y$ play
important roles. Another important thing is that we can write
$Z_4(A,A)$ in two ways such that the third coordinate of one of
these (i.e. $(abx+1)/aA(x)$) is mapped to the third coordinate of
$Z_2(A,A)$ (i.e. $(x-1)/x$) under $\rho_x$ and vice versa. Hence
$$Z(A,A)=\sum_{i=1}^4 Z_i(A,A).$$
On the other hand, we can easily see that
\begin{equation}\label{Z3fB}
\rho_{x,y} Z_1(A,A)= Z_3\Bigl(\frac{f}{B},\frac{f}{B}\Bigr):=
\Bigl[\frac{(b-1)x}{B(x)},\frac{(1-b)y}{B(y)},\frac{abx+1}{
aA(x)},
    \frac{y-x}{yB(x)},l(y)\Bigr].
\end{equation}
Therefore we have the following simple expression of $\{k(c)\}$ by
\eqref{AA}:
\begin{equation}\label{kc}
\{k(c)\}=\sum_{i=1}^4 Z_i(A,A)
=Z_3(A,A)+Z_3\Bigl(\frac{f}{B},\frac{f}{B}\Bigr)
+\rho_x  Z_2(A,A)+\rho_y Z_4(A,A)
\end{equation}
where
\begin{align}
\rho_x Z_2(A,A)=&\Bigl[\frac {(b-1)x}{B(x)},\frac{(b-1)A(y)}{\mu},
        \frac{abx+1}{aA(x)},\frac{y-x}{A(y)},l(y)\Bigr],
                \label{rhox}\\
\rho_y  Z_4(A,A)=&\Bigl[\frac{(b-1)A(x)}{\mu},
\frac{(b-1)y}{B(y)},
    \frac{abx+1}{aA(x)},\frac{\mu(x-y)}{A(y)B(x)},l(y)\Bigr].
        \label{rhoy}
\end{align}

\medskip
\noindent
{\bf Step} (4). Decomposition of
$\rho_x Z_2(A,A)+\rho_y Z_4(A,A)=X_1-X_2$.

\medskip

Let $f_1(x)=(b-1)x/B(x)$, $f_3(x)=(abx+1)/aA(x)$, $f_2=gh$ where
$g(x)=A(x)/(-\mu x)$ and $h(x)=(1-b)x$. Then we can apply
Lemma~\ref{keylem}(ii) and easily get
\begin{equation}\label{rho24}
\rho_x Z_2(A,A)+\rho_x Z_4(A,A) = X_1-X_2
\end{equation}
where
$$\aligned
X_1=&\Bigl[\frac{(b-1)x}{B(x)} ,\frac{A(y)}{-\mu y},
    \frac{abx+1}{aA(x)},\frac{y-x}{A(y)},l(y)\Bigr]\\
+&\Bigl[\frac{A(y)}{-\mu x},\frac{(b-1)y}{B(y)},
\frac{abx+1}{aA(x)},
    \frac{\mu(x-y)}{A(y)B(x)},l(y)\Bigr],\\
X_2=&\Bigl[\frac{B(x)}{(b-1)x},(1-b)y,\frac{abx+1}{aA(x)},
    \frac{y-x}{A(y)},l(y)\Bigr]\\
+&\Bigl[(1-b)x,\frac{B(y)}{(b-1)y}, \frac{abx+1}{aA(x)},
    \frac{\mu(x-y)}{A(y)B(x)},l(y)\Bigr].
\endaligned$$
We now apply Lemma~\ref{keylem}(ii) to $X_2$ with $g=-1$ and
$h(x)=(b-1)x$ to get
\begin{align*}
X_2=&\Bigl[\frac{B(x)}{(b-1)x},(b-1)y,\frac{abx+1}{aA(x)},
    \frac{y-x}{A(y)},l(y)\Bigr]\\
+&\Bigl[(b-1)x,\frac{B(y)}{(b-1)y}, \frac{abx+1}{aA(x)},
    \frac{\mu(x-y)}{A(y)B(x)},l(y)\Bigr].
\end{align*}

\medskip
\noindent
{\bf Step} (5). Computation of  $X_1$.

\medskip

Set
$$\tZ(f_1,f_2)=\Bigl[f_1,f_2,\frac{abx+1}{aA(x)},
    \frac{\mu(x-y)}{A(y)B(x)},l(y)\Bigr].$$
Throwing away the appropriate admissible and negligible cycle we have
$$X_1=\tZ\Bigl(\frac{(b-1)f}{B},\frac{A}{-\mu f}\Bigr)
    +\tZ\Bigl(\frac{A}{-\mu f},\frac{(b-1)f}{B}\Bigr).$$
Then further disregarding some admissible and negligible cycles we get
$$Z_3(F,F)=\tZ(F,F) \text{ for }F=\frac{A}{f},\
    \frac{f}{B},\ \frac{A}{B}$$
where $Z_3(f/B,f/B)$ is defined by
\eqref{Z3fB},
\begin{equation}\label{Z3Af}
Z_3\Bigl(\frac{A}{f},\frac{A}{f}\Bigr)
 :=\Bigl[\frac{A(x)}{-\mu x},\frac{A(y)}{-\mu y},
    \frac{abx+1}{aA(x)}, \frac{y-x}{A(y)}, l(y)\Bigr]
\end{equation}
and
\begin{equation}\label{Z3AB}
Z_3\Bigl(\frac{A}{B},\frac{A}{B}\Bigr):=\Bigl[\frac{A(x)}{B(x)},
    \frac{A(y)}{B(y)},\frac{abx+1}{aA(x)},
\frac{-\mu y}{A(y)B(x)}\Bigl(1-\frac{x}{y}\Bigr),l(y)\Bigr].
\end{equation}
Here we removed the coefficient $b-1$ in front of the $A/B$ and
$-1/\mu$ in front of $A/f$ in $Z_3$ by Lemma \ref{fgh}(ii). Then
we can take $f_1=f_2=(b-1)A/B$, $g=A/f$ and $h=(b-1)f/B$ in Lemma
\ref{fgh}(ii)(c) and get
\begin{align}
X_1=&\tZ\Bigl(\frac{A}{B},\frac{A}{B}\Bigr)
-\tZ\Bigl(\frac{A}{f},\frac{A}{f}\Bigr)
-\tZ\Bigl(\frac{f}{B},\frac{f}{B}\Bigr)\nonumber \\
=&Z_3\Bigl(\frac{A}{B},\frac{A}{B}\Bigr)
-Z_3\Bigl(\frac{A}{f},\frac{A}{f}\Bigr)
-Z_3\Bigl(\frac{f}{B},\frac{f}{B}\Bigr).\label{X1}
\end{align}

\medskip
\noindent
{\bf Step} (6). Decomposition of $X_2=Y_1+Y_2+Y_3+Y_4$.

\medskip
Put
$$v(x)=\frac{abx+1}{aA(x)},
\quad  l_1(y)=1-\frac{y}{c}, \quad l_2(y)=\frac{y_2-y}{y_2B(y)},
$$ which satisfies
$$l_1(y)l_2(y)=l(y)=1-\frac{k(c)}{k(y)},\quad l_1(0)=l_2(0)=1.$$
Then it follows from Lemma~\ref{keylem}(i) that
\begin{equation}\label{X2}
X_2=Y_1+Y_2+Y_3+Y_4
\end{equation}
where all of the cycles
\begin{align*}
Y_1=&\Bigl[\frac{B(x)}{(b-1)x}, (b-1)y,\frac{abx+1}{aA(x)},
    \frac{y-x}{A(y)},l_1(y)\Bigr],\\
Y_2=&\Bigl[(b-1)x,\frac{B(y)}{(b-1)y},\frac{abx+1}{aA(x)},
    \frac{\mu(x-y)}{A(y)B(x)},l_1(y)\Bigr],\\
Y_3=&\Bigl[\frac{B(x)}{(b-1)x}, (b-1)y,\frac{abx+1}{aA(x)},
    \frac{y-x}{A(y)},l_2(y)\Bigr],\\
Y_4=&\Bigl[(b-1)x,\frac{B(y)}{(b-1)y},\frac{abx+1}{aA(x)},
    \frac{\mu(x-y)}{A(y)B(x)},l_2(y)\Bigr]
\end{align*}
are admissible. This breakup is the key step in the whole paper.

\medskip
\noindent
{\bf Step} (7). Computation of $Y_1+Y_2$.

\medskip

To ease the reading of the proof in this section we first
summarize our approach here. We very much like to be able to use
Lemma~\ref{fgh}(ii) but unfortunately the terms $f_3$, $f_4$ and
$f_5$ cannot be fixed for all the terms because we have to stay
inside the ``admissible world''. Nevertheless, luckily enough for
us, most of the cycles we are going to use have more than one
``realization'' so that we can apply Lemma~\ref{fgh} and
Lemma~\ref{keylem} to obtain the desired results. Corollary
\ref{key} will be crucial to our computation.

We begin by setting
$$\ga=\frac{bc-c}{bc-c+1}, \qquad \gd= \frac{1}{b},$$
and
$${\setlength\arraycolsep{1pt}
\begin{array}{rlrlrl}
v(x)&={\displaystyle \frac{abx+1}{aA(x)}}, & \quad
g(x)&={\displaystyle \frac{B(x)}{(b-1)x}}, &\quad h(x)&=(b-1)x,
\\
\quad p_4(x,y)&={\displaystyle \frac{\mu(x-y)}{A(y)B(x)}}, & \quad
q_4(x,y)&={\displaystyle \frac{y-x}{A(y)}}, & \quad \quad
s_4(x,y)&={\displaystyle \frac{(b-1)(y-x)}{B(y)}},
\\
r_4(x,y)&={\displaystyle \frac{(b-1)(y-x)}{xB(y)}}, & \quad
w_4(x,y)&={\displaystyle \frac{y-x}{B(x)(y-1)}}.
\end{array}}$$
such that $\ga l_1\big(1/(1-b)\big)=\gd v(\infty)=1$. By Lemma
3.1(ii)(1) we get
\begin{align*}
2[gh, gh,\gd v, q_4,\ga l_1]
=&[gh, gh,\gd v, q_4,\ga l_1]+[gh, gh,\gd v, s_4,\ga l_1]\\
 =&[g, gh,\gd v, q_4, \ga l_1]+[h, gh,\gd v, q_4, \ga l_1]\\
    \ +&[gh, g,\gd v, s_4, \ga l_1]+[gh, h,\gd v, s_4, \ga l_1]
\end{align*}
are all admissible. Then repeatedly applying Lemma~\ref{fgh} and
Lemma~\ref{keylem} we have
\begin{align*}
 \ &[g, gh,\gd v, q_4, \ga l_1]+[gh, g,\gd v, s_4, \ga l_1]\\
 =&[g,gh, \gd v, q_4, \ga l_1]+[gh, g, \gd v, r_4, \ga l_1]\\
 =&[g,gh, v, q_4, \ga l_1]+[gh, g,  v, r_4, \ga l_1]\\
 =&[g,gh, v, q_4, \ga l_1]+[gh, g,  v, w_4, \ga l_1]\\
 =&[g, gh,v, q_4,  l_1]+[gh, g, v, w_4, l_1]\\
 =&[g, gh,v, q_4,  l_1]+[gh, g, v, p_4, l_1]\\
 =&[g,h,v, q_4, l_1]+[h, g, v, p_4, l_1]
 +[g, g, v,q_4,l_1]+[g,g, v, p_4,l_1]\\
 =&[g,h,v, q_4, l_1]+[h, g, v, p_4, l_1]
 +2[g, g, v, p_4,l_1].
\end{align*}
Again by applying Lemma~\ref{fgh} and Lemma~\ref{keylem} we get
\begin{align*}
 \ &[h, gh,\gd v,q_4,\ga l_1]+[gh, h,\gd v, s_4,\ga l_1]\\
 =&[h, gh,\gd v,q_4, l_1]+[gh, h,\gd v, s_4, l_1]\\
 =&[h, gh,\gd v,q_4, l_1]+[gh, h,\gd v, q_4, l_1]\\
 =&[h, g, \gd v,q_4, l_1]+[g, h,\gd v, q_4, l_1]
 +2[h, h,\gd v,q_4, l_1]\\
 =&[h, g, v, q_4, l_1]+[g, h, v, q_4, l_1]
 +2[h, h,\gd v,q_4, l_1]\\
 =&[h, g, v, p_4, l_1]+[g, h, v, q_4, l_1]
 +2[h, h,\gd v,q_4, l_1].
\end{align*}
Therefore
\begin{align}\nonumber
Y_1+Y_2=&[h, g, v, p_4, l_1]+[g, h, v, q_4, l_1]\\
=&[gh, gh,\gd v, q_4, \ga l_1]- [g, g, v, p_4,l_1]
    -[h, h,\gd v,q_4, l_1].\label{Y12}
\end{align}

\medskip
\noindent {\bf Step} (8). Computation of $Y_3+Y_4$.
\medskip

We could use a similar process as in Step (7) to do the computation.
But we can get around this by the following argument. Define
the substitutions
$$\gs_{x,y}: (x,y)\lms \Bigl(\frac{-x}{B(x)},\frac{-y}{B(y)}\Bigr),
\quad \tau_{a,c}: (a,c)\lms \Bigl(\frac{ab-b+1}{b(a-1)},
 \frac{ca-a+1}{ab-b+1}\Bigr).$$
Let
\begin{align*}
Y_3'=&\Bigl[\frac{(1-b)x}{B(x)},
 (1-b)y,\frac{abx+1}{abA(x)},\frac{y-x}{A(y)},l_2(y)\Bigr],\\
Y_4'=&\Bigl[(1-b)x,\frac{(1-b)y}{B(y)},\frac{abx+1}{abA(x)},
    \frac{\mu(x-y)}{A(y)B(x)},l_2(y)\Bigr].
\end{align*}
Then an easy computation shows that
\begin{align*}
-Y_4' =&\tau_{a,c}\gs_{x,y} Y_1 =
\Bigl[\frac{1}{(1-b)x},\frac{(1-b)y}{B(y)},\frac{abA(x)}{abx+1},
 \frac{(ab-b+1)(y-x)}{(aby+1)B(x)},l_2(y)\Bigr],\\
-Y_3' =&\tau_{a,c}\gs_{x,y}Y_2 = \Bigl[\frac{(1-b)x}{B(x)},
 \frac{1}{(1-b)y},\frac{abA(x)}{abx+1},
 \frac{ab(y-x)}{aby+1},l_2(y)\Bigr]
\end{align*}
Hence by first splitting off the $-1$ in front of the first two
coordinates of $Y_3'$ and $Y_4'$ respectively, then removing some
other admissible and negligible cycles we find
\begin{equation}\label{Y34}
Y_3+Y_4=Y'_3+Y'_4=-\tau_{a,c}(Y_1+Y_2).
\end{equation}

\medskip
\noindent {\bf Step} (9). Final decomposition of $\{k(c)\}$ into
$T_i(F)$'s.

\medskip
Putting \eqref{kc} to \eqref{Y34} together we see that
\begin{align}
\{k(c)\}=&Z_3(A,A) + Z_3\Bigl(\frac{A}{B},\frac{A}{B}\Bigr)
-Z_3\Bigl(\frac{A}{f},\frac{A}{f}\Bigr)\notag  \\
+&  (1-\tau_{a,c})\big([g, g, v, p_4, l_1] +[h, h,\gd v, q_4, l_1]
-[gh, gh,\gd v, q_4, \ga l_1]\big).  \label{needsim}
\end{align}
We will first simplify the terms in the above expression. Set
$$\eps_1(f)=\eps_2(f)=1,\quad
\eps_1(A)=\frac{ca}{ca-a+1},\quad
\eps_2(A)=\frac{ca-a+1}{ca}.
$$
Define the admissible cycles
$$T_i(F)=
\begin{cases}
{\displaystyle \Bigl[\frac{A(x)}{x},\frac{A(y)}{y},
    \frac{(1-a)(x-1)}{x}, 1-\frac{x}{y}, l_i(y)\Bigr]} &
    \quad \text{if }F=\frac{A}{f},\ i=1,2, \\
{\displaystyle \Bigl[F(x),F(y), \frac{x-1}{x}, \frac{y-x}{F(y)},
    \eps_i(F) l_i(y)\Bigr]} &
    \quad \text{if }F=f, A,\ i=1,2, \\
{\displaystyle \Bigl[F(x),F(y), \frac{abx+1}{abA(x)},
    \frac{y-x}{F(y)}, \eps_{i-2}(F)l_{i-2}(y)\Bigr]} &  \quad
    \begin{cases}
        \text{if }& F=A,\ i=3,4,\\
        \text{if }& F=f,\ i=3,
    \end{cases} \\
{\displaystyle \Bigl[B(x),B(y),\frac{abx+1}{abA(x)},
    \frac{y-x}{A(y)},\ga l_1(y)\Bigr]}
        &\quad \text{if }F= B, i=1,  .
\end{cases}$$

\noindent{\bf Claim}. We have
$$\{k(c)\}=\sum_{i=1}^3 T_i(f)+\sum_{i=2}^4 T_i(A)
-\sum_{i=1,2} T_i\Bigl(\frac{A}{f}\Bigr)
-T_1(B)-\tau_{a,c}\big(T_3(f)+T_2(A)-T_1(B)\big).$$

\noindent{\em Proof of the Claim. }
First by the involution we quickly find
$$Z_3(A,A)=T_3(A)+T_4(A), \quad
Z_3\Bigl(\frac{A}{B},\frac{A}{B}\Bigr)=T_1(f)+T_2(f). $$ For
$Z_3(A/f,A/f)$ defined by \eqref{Z3Af} we can remove coefficient
$-1/\mu$ from the first two coordinates then apply $\rho_{x,y}$ to
get
$$\aligned
Z_3\Bigl(\frac{A}{f},\frac{A}{f}\Bigr)
=&\Bigl[\frac{A(x)}{x},\frac{A(y)}{y},
    \frac{x-1}{x}, \frac{y-x}{yB(x)},l(y)\Bigr]\\
=&\Bigl[\frac{A(x)}{x},\frac{A(y)}{y},
    \frac{x-1}{x}, 1-\frac{x}{y},l(y)\Bigr]\\
=&\Bigl[\frac{A(x)}{x},\frac{A(y)}{y},
    \frac{(1-a)(x-1)}{x}, \frac{y-x}{yB(x)},l(y)\Bigr]\\
=&T_1\Bigl(\frac{A}{f}\Bigr)+T_2\Bigl(\frac{A}{f}\Bigr).
\endaligned$$
Here we have sequentially added two admissible and negligible cycles.

For the first term in \eqref{needsim}, from \eqref{rhoxys} and
$l_1(y)\os{\rho_{x,y}}{\llra} \eps_2(A)l_2(y)$ we find
$$\rho_{x,y} [g, g, v, p_4, l_1]
=\Bigl[A(x),A(y), \frac{x-1}{x}, 1-\frac{x}{y},
    \eps_2(A)l_2(y)\Bigr] =T_2(A).$$
For the second term in \eqref{needsim}, we can first remove the
coefficient $b-1$ in $h$ and then delete one admissible and
negligible cycle to get $[h, h,\gd v, q_4, l_1]=T_3(f).$ This
completes the proof of our claim.

\medskip
\noindent {\bf Step} (10). Final computation of $\{k(c)\}$.

\medskip

Let's compute each $T_i(F)$ separately. Throughout this
computation we will repeatedly invoke Lemma \ref{stuv} without
explicitly stating it.

\medskip \noindent \boxed{\text{$F=f$}} \par\medskip By definition
$$\aligned
T_1(f)=&\{c\}, \quad
T_2(f)=\Bigl\{\frac{-a(bc-c+1)}{ca-a+1}\Bigr\}-\{1-b\}, \\
T_3(f)=&\{-abc\}-\Bigl\{1-\frac{ca-a+1}{ca} \Bigr\},  \\
\tau_{a,c}T_3(f)=&\Bigl\{1-\frac{ca}{ca-a+1}\Bigr\}-\{ca-a+1\}.
\endaligned $$

\medskip \noindent \boxed{\text{$F=A$}}

\medskip

Using $(x,y)\mapsto (x+(a-1)/a,y+(a-1)/a)$ we find
$$\aligned
T_2(A)=&\Bigl\{\frac{-c(ab-b+1)}{bc-c+1}\Bigr\}
-\Bigl\{1-\frac{ab}{ab-b+1}\Bigr\} \\
\ &\hskip15ex - \Bigl\{1-\frac{ca-a+1}{c(ab-b+1)} \Bigr \}+
\Bigl\{1-\frac{a}{ab-b+1}\Bigr\},\\
\tau_{a,c}T_2(A)=&
\Bigl\{1-\frac{c(ab-b+1)}{ca-a+1}\Bigr\}
-\Bigl\{1-\frac{ab-b+1}{a}\Bigr\} \\
\ &\hskip15ex -\Bigl\{\frac{bc-c+1}{b(ca-a+1)}\Bigr\}
+\Bigl\{1-\frac{ab-b+1}{ab}\Bigr\},
\endaligned$$
$$T_3(A)=\Bigl\{\frac{-b(ca-a+1)}{ab-b+1}\Bigr\},\quad
T_4(A)=\Bigl\{\frac{bc-c+1}{bc}\Bigr\}-\Bigl\{1-\frac{1}{b}\Bigr\}.
$$

\medskip \noindent \boxed{\text{$F=A/f$}} \par\medskip
Using substitution $(x,y)\mapsto ((1-a)/(ax-a),(1-a)/(ay-a) )$  we get
$$\aligned
T_1\Bigl(\frac{A}{f} \Bigr)
=&\Bigl\{\frac{ca-a+1}{c}\Bigr\}-\{a\},\\
T_2\Bigl(\frac{A}{f} \Bigr)
=&\Bigl\{\frac{ca(ab-b+1)}{ca-a+1} \Bigr\}-\{ab-b+1\}.
\endaligned$$

\medskip \noindent \boxed{\text{$F=B$}} \par\medskip
By definition and using substitution $(x,y)\mapsto
((x-1)/(b-1),(y-1)/(b-1)$ we get
$$\aligned
T_1(B)=& \Bigl\{\frac{ab-b+1}{ab(bc-c+1)}\Bigr\}
    -\Bigl\{\frac{ab-b+1}{a(bc-c+1)}\Bigr\},\\
\tau_{a,c}T_1(B)=& \{bc-c+1\} - \Bigl\{\frac{bc-c+1}{b}\Bigr\}.
\endaligned$$

Putting the above together we now complete the proof the theorem
in the case that none of the terms in Goncharov's relations is
equal to $\{0\}$ or $\{1\}$.
\end{proof}

\end{document}